# Algorithm 916: computing the Faddeyeva and Voigt functions

MOFREH R. ZAGHLOUL, United Arab Emirates University[1]

AHMED N. ALI, United Arab Emirates University

We present a MATLAB function for the numerical evaluation of the Faddeyeva function $w(z)$. The function is based on a newly developed accurate algorithm. In addition to its higher accuracy, the software provides a flexible accuracy vs efficiency trade-off through a controlling parameter that may be used to reduce accuracy and computational time and vice versa. Verification of the flexibility, reliability and superior accuracy of the algorithm is provided through comparison with standard algorithms available in other libraries and software packages.

---

Author's addresses: M. Zaghloul and A. Ali, Department of Physics, College of Sciences, United Arab Emirates University, Al-Ain, 17551, UAE.

## 1. INTRODUCTION

The value and significance of the scaled complementary error function for complex variables, also known as the Faddeyeva function or the plasma dispersion function, is well recognized in the literature for its applications in several fields of physics [Armstrong 1967; Gautschi 1967; 1970; Hui et al. 1978; Humlíček 1982; Dominguez et al. 1987; Poppe et al. 1990; Lether et al. 1991; Schreier 1992; Shippony et al. 1993; Weideman 1994; Wells 1999; Luque et al. 2005; Letchworth et al. 2007; Abrarov et al. 2010]. Plasma spectroscopy, nuclear physics, radiative heat transfer, and nuclear magnetic resonance are a few examples of the fields for which efficient and accurate evaluation of this function is required. Some of these applications require a small number of evaluations of the function where accuracy is more important than the computational time while other applications require enormous numbers of function evaluations, which imposes tight restrictions on the computational time. Accordingly, computational accuracy and computational time are issues of interest that should be critically addressed and investigated in developing any successful algorithm for the computation of this function.

Motivated by its practical importance and a lack of closed form expressions for the calculation of the Faddeyeva function, numerical evaluation of the function has been the focus of research over many decades [Armstrong 1967; Gautschi 1967; 1970; Hui et al. 1978; Humlíček 1982; Dominguez et al. 1987; Poppe et al. 1990; Lether et al. 1991; Schreier 1992; Shippony et al. 1993; Weideman 1994; Wells 1999; Luque et al. 2005; Letchworth et al. 2007; Abrarov et al. 2010; Zaghloul 2007]. As a result, a wide variety of algorithms for the calculation of this function have been developed and presented in the literature. However, as it is shown in this study, most of these algorithms lose accuracy in some regions of the computational domain.

We introduce a new algorithm for the calculation of the Faddeyeva function which provides flexibility, reliability and superior accuracy. In section 2 we present the definition of the function and briefly summarize some relevant fundamental mathematical relations. Then, in section 3, we establish the analytical basis of the algorithm while numerical analysis and computational details are discussed in section 4. A short description of the Matlab function is given in section 5. Verification of the algorithm and comparisons with other competitive codes in the literature are provided in section 6.

## 2. DEFINITION AND FUNDAMENTAL MATHEMATICAL RELATIONS

For a complex variable $z=x+iy$, the Faddeyeva (plasma dispersion) function, $w(z)$, the real Voigt function, $V(x,y)$, the imaginary Voigt function, $L(x,y)$, the complex error function, erf($z$), the imaginary error function, erfi($z$) and the Dawson's integral $F(z)$ are all closely related to each other. One can summarize these relations as

$$\begin{aligned} w(z) &= e^{-z^2}(1-\text{erf}(-iz)) = e^{-z^2}\text{erfc}(-iz) \\ &= e^{-z^2}(1+\text{erf}(iz)) \\ &= e^{-z^2}(1+i\,\text{erfi}(z)) \\ &= e^{-z^2} + \frac{2i}{\sqrt{\pi}}F(z) \\ &= V(x,y) + i\,L(x,y) \qquad \text{for } y>0 \end{aligned} \qquad (1)$$

In the above relations, $i = \sqrt{-1}$, erfc(z) is the complementary error function, erfi(z) is the imaginary error function which is related to the error function by erfi(z) = -i erf(iz).

As can be seen from the last line in (1), the real and imaginary Voigt functions are just the real and imaginary parts of the Faddeyeva function for $y > 0$, respectively. The evaluation of all of these functions can, therefore, be performed through the error function.

The error function of a complex variable $z$ can be regarded as a line integral in the complex plane given by

$$\mathrm{erf}(z) = \frac{2}{\sqrt{\pi}} \int_0^z e^{-t^2}\, dt \tag{2}$$

Different paths can be taken to perform this line integral in the complex plane. For example, one may choose a linear path between the initial point (origin) and the final point $z$ which gives an expression for the complex error function of the form

$$\mathrm{erf}(z) = \frac{2z}{\sqrt{\pi}} \int_0^1 e^{-z^2 t^2}\, dt \tag{3}$$

An alternative path can be followed through the line segments $[0, 0 \rightarrow iy]$ and $[0 \rightarrow x, iy]$ which gives, for the complex error function the expression,

$$\mathrm{erf}(z) = \frac{2 e^{y^2}}{\sqrt{\pi}} \int_0^x e^{-t^2} \cos(2yt)\, dt + i \left[ \frac{2}{\sqrt{\pi}} \int_0^y e^{t^2}\, dt - \frac{2 e^{y^2}}{\sqrt{\pi}} \int_0^x e^{-t^2} \sin(2yt)\, dt \right] \tag{4}$$

In addition to the above paths, another simple and useful path from the initial to the final points can be taken through the line segments $[0 \rightarrow x, iy=0]$ and $[x, 0 \rightarrow iy]$ which results in

$$\mathrm{erf}(z) = \frac{2}{\sqrt{\pi}} \int_0^x e^{-t^2}\, dt + e^{-x^2} \frac{2}{\sqrt{\pi}} \int_0^y e^{t^2} \sin(2xt)\, dt + i\, e^{-x^2} \frac{2}{\sqrt{\pi}} \int_0^y e^{t^2} \cos(2xt)\, dt \tag{5}$$

The first term on the right hand side of (5) is the definition of the error function of the real variable $x$. Equation (5) is the basis of the present algorithm for the calculation of the Faddeyeva function as shown below.

For some limiting values of the parameters $x$ and $y$, analytical formulae do exist for the real and imaginary parts of the Faddeyeva function

$$V(x \rightarrow 0, y) \rightarrow \mathrm{erfcx}(y)$$
$$L(x \rightarrow 0, y) \rightarrow 0$$
$$V\!\left((x^2 + y^2) \rightarrow \infty\right) \rightarrow (1/\sqrt{\pi})\left[y/(x^2 + y^2)\right] \tag{6}$$
$$L\!\left((x^2 + y^2) \rightarrow \infty\right) \rightarrow (1/\sqrt{\pi})\left[x/(x^2 + y^2)\right]$$
$$V(\pm x, y \rightarrow 0) \rightarrow e^{-x^2}$$

where $\mathrm{erfcx}(y) = e^{y^2} \mathrm{erfc}(y)$ is the scaled complementary error function of the real argument $y$. When $y \rightarrow 0$, the imaginary part of the Faddeyeva function cannot be expressed as simply, however, it can be expressed in terms of Dawson's integral of $x$ (the real part of $z$) where

$$L(\pm x, y \rightarrow 0) \rightarrow 2/\sqrt{\pi}\, F(x) \tag{7}$$

and $F(x)$ has well reported asymptotic expressions for limiting values of $x$ [Armstrong 1967].

Following Salzer 1951 [Salzer 1951], we can write

$$e^{t^2} \cong \frac{a}{\sqrt{\pi}} (1 + 2\sum_{n=1}^{\infty} e^{-a^2 n^2} \cosh(2ant)) \pm E\, e^{t^2} \qquad\qquad a \leq 1 \qquad (8)$$

where the relative error $E = 2\sum_{n=1}^{\infty} e^{-n^2\pi^2/a^2} \cos(2n\pi t/a)$ is of the order of $E \sim 2e^{-\pi^2/a^2}$.

Table 1 presents some representative values of the relative error $E$ corresponding to some values of the parameter $a$.

Table 1: Relative error, $E$, of the representation of $e^{t^2}$ given by Eq. (8) as a function of the parameter $a$.

| $a$ | 1.0 | 0.9 | 0.8 | 0.7 | 0.6 | 0.5 | 0.4 | 0.36 |
|---|---|---|---|---|---|---|---|---|
| $E$ | $1.03\times10^{-4}$ | $1.02\times10^{-5}$ | $4.01\times10^{-7}$ | $3.58\times10^{-9}$ | $2.48\times10^{-12}$ | $1.43\times10^{-17}$ | $3.25\times10^{-27}$ | $1.69\times10^{-33}$ |

## 3- ANALYTICAL BASIS OF THE ALGORITHM

Replacing $e^{t^2}$ in (5) by its representation given in (8), we get, for the real and imaginary parts of erf($z$), the following expressions

$$\mathrm{Re}[\mathrm{erf}(z)] = \mathrm{erf}(x) + \frac{a}{\pi x} e^{-x^2}(1 - \cos(2xy)) + \frac{8a}{\pi} e^{-x^2}$$
$$\times \sum_{n=1}^{\infty} \frac{e^{-a^2 n^2}}{4a^2 n^2 + 4x^2}\left[x - x\cosh(2any)\cos(2xy) + an\sinh(2any)\sin(2xy)\right] \qquad (9)$$

$$\mathrm{Im}[\mathrm{erf}(z)] = \frac{4a}{\pi} e^{-x^2}(4x)^{-1}\sin(2xy) + \frac{8a}{\pi} e^{-x^2}$$
$$\times \sum_{n=1}^{\infty} \frac{e^{-a^2 n^2}}{4a^2 n^2 + 4x^2}\left[x\cosh(2any)\sin(2xy) + an\sinh(2any)\cos(2xy)\right] \qquad (10)$$

The error in the expressions in (9) and (10) for the real and imaginary parts of the complex error function is controllable through the parameter $a$ as indicated in Table 1. Both expressions, however, reduce to the formulae given by Salzer [Salzer 1951] and Abramowitz and Stegun [Abramowitz et al. 1972] with a relative error less than the floating point relative accuracy, $\varepsilon$, on a 16 digit computational platform, by setting $a$ equal to $1/2$. An expression for the Faddeyeva function, $w(z) = \mathrm{Re}[w(z)] + i\,\mathrm{Im}[w(z)] = e^{-z^2}\mathrm{erfc}(-iz)$, can be obtained by substituting (9) and (10) into (1). The resulting expressions for the real and imaginary parts of the Faddeyeva function are then given by

$$\mathrm{Re}[w(z)] = e^{-x^2}\mathrm{erfcx}(y)\cos(2xy) + \frac{2a\,x\sin(xy)}{\pi}e^{-x^2}[\sin(xy)/(xy)]$$

$$+ \frac{8a\,y}{\pi}e^{-x^2}\sum_{n=1}^{\infty}\frac{e^{-a^2n^2}}{4a^2n^2+4y^2}[\cosh(2anx)-\cos(2xy)]$$

(11)

$$\mathrm{Im}[w(z)] = -e^{-x^2}\mathrm{erfcx}(y)\sin(2xy) + \frac{2a\,x}{\pi}e^{-x^2}[\sin(2xy)/(2xy)]$$

$$+ \frac{8a}{\pi}e^{-x^2}\sum_{n=1}^{\infty}\frac{e^{-a^2n^2}}{4a^2n^2+4y^2}[y\sin(2xy)+a\,n\sinh(2anx)]$$

(12)

Some comments on the above expressions may be useful at this stage
1- The series in (11) and (12) have an infinite number of terms and need to be truncated for practical use. The effect of such a truncation on the accuracy of the computations is predictable and can be controlled for a converging series,
2- While the exponential factors in the series are decaying with the index, *n*, they are multiplied by growing hyperbolic functions and one cannot directly determine where to truncate the infinite sums for practical use,
3- There is a limitation on the evaluation of the hyperbolic functions (related to the largest positive floating point number, $R_{max}$, available on the computational platform). Since the argument of these hyperbolic functions is *(2anx)*, this will impose a strict restriction on the number of terms to be included in the sums for a given value of *x* with possible catastrophic consequences on both accuracy and reliability,
4- In writing the above expressions, the scaled complementary error function of a real variable is used to reduce rounding error associated with the term (1-erf(*y*)), however, special care is needed for the evaluation of the erfcx(*y*) function to avoid overflow problems illustrated in [Zaghloul 2007]. At present, many software packages have well-behaved algorithms for computing the scaled complementary error function of a real variable, erfcx, and algorithms for accurate and efficient computation of this function are available in the literature,
5- The evaluation of the quantity exp(-$x^2$) common to all of the above terms, can suffer underflow problems for large values of *x*. These problems can be avoided by combining this quantity with other large quantities wherever possible.

To overcome the above-stated concerns and restrictions, the expressions in (11) and (12) are rewritten in the forms

$$\mathrm{Re}[w(z)] = e^{-x^2}\mathrm{erfcx}(y)\cos(2xy) + \frac{2a\,x\sin(xy)}{\pi}e^{-x^2}[\sin(xy)/(xy)]$$

$$+ \frac{2a}{\pi}\left\{-y\cos(2xy)\Sigma_1 + \frac{y}{2}\Sigma_2 + \frac{y}{2}\Sigma_3\right\}$$

(13)

and

$$\text{Im}[w(z)] = -e^{-x^2}\text{erfcx}(y)\sin(2xy) + \frac{2ax}{\pi}e^{-x^2}[\sin(2xy)/(2xy)]$$

$$+ \frac{2a}{\pi}\left\{y\sin(2xy)\,\Sigma_1 - \frac{1}{2}\Sigma_4 + \frac{1}{2}\Sigma_5\right\} \tag{14}$$

where

$$\Sigma_1 = \sum_{n=1}^{\infty}\left(\frac{1}{a^2n^2+y^2}\right)e^{-(a^2n^2+x^2)} \tag{15}$$

$$\Sigma_2 = \sum_{n=1}^{\infty}\left(\frac{1}{a^2n^2+y^2}\right)e^{-(an+x)^2} \tag{16}$$

$$\Sigma_3 = \sum_{n=1}^{\infty}\left(\frac{1}{a^2n^2+y^2}\right)e^{-(an-x)^2} \tag{17}$$

$$\Sigma_4 = \sum_{n=1}^{\infty}\left(\frac{an}{a^2n^2+y^2}\right)e^{-(an+x)^2} \tag{18}$$

$$\Sigma_5 = \sum_{n=1}^{\infty}\left(\frac{an}{a^2n^2+y^2}\right)e^{-(an-x)^2} \tag{19}$$

The convergence of the series (15)-(19) can be verified by applying the simple *ratio-test* [Boas 2006]. In addition, in all of the above expressions (15)-(19), the pre-exponential factors (fractions in brackets) of the arguments of the summations assume values less than 1 for $n>1/a$. For $a=1/2$ (which is sufficient to express $e^{t^2}$ by the expression given in (8) to machine accuracy on a 16-digit computational platform, see Table 1) the pre-exponential factors will always assume values less than or equal to one for $n\geq 2$.

For values of $n\geq 2$, the terms in the series $\Sigma_1$, $\Sigma_2$ and $\Sigma_4$ decrease monotonically with $n$ while for relatively large values of $x$ the terms $\Sigma_3$ and $\Sigma_5$ increase with $n$ from one to a certain limit when they start decaying monotonically with $n$. The fact that all exponential factors in the above summations eventually decay with increasing $n$ provides us with the possibility of obtaining a truncated series of practical use and computational efficiency as shown in the next section. It is understood, however, that all of the above summations are to be performed using a single loop for computational efficiency. The loop specifications can be determined once a cutoff scheme for the series in (15)-(19) is established. Moreover the computation of the exponentials in these series can be reduced to the computation of a single exponential and products using a single computational loop. This reduces the dependence on using the intrinsic function to calculate these exponentials within the loop and can save a significant amount of computational time.

As the real part of the Faddeyeva function is even in $x$ and its imaginary part is odd in $x$, we need only consider the right half of the complex plane ($x\geq 0$) since the even/odd properties of the real and imaginary parts of the function can be used to find the corresponding values in the left half of the plane. In addition, the values of $w(z)$ in the lower half of the complex plane can be obtained from values in the upper half using the relationship

$$w(-z) = 2e^{-z^2} - w(z) \tag{20}$$

which is equivalent to the symmetry relations

$$V(x,y) = 2e^{-x^2+y^2}\cos(2xy) - V(-x,-y)$$
$$L(x,y) = -2e^{-x^2+y^2}\sin(2xy) - L(-x,-y)$$
(20')

Thus we do not lose generality by considering the evaluation of the function in just the first quadrant, however, we note that (20) requires the subtraction of two values which leads to a loss of accuracy in the lower half.

The practical usefulness of the calculation of the partial derivatives of the real and imaginary parts of the Faddeyeva function has been pointed out by many authors [10,13,15]. Once the Faddeyeva function has been calculated accurately one can calculate these partial derivatives with relative simplicity using the expressions

$$\frac{\partial V(x,y)}{\partial x} = -2\,\mathrm{Re}[z\,w(z)] = 2[y\,L(x,y) - x\,V(x,y)] \tag{21}$$

$$\frac{\partial V(x,y)}{\partial y} = 2\,\mathrm{Im}[z\,w(z)] - \frac{2}{\sqrt{\pi}} = 2[x\,L(x,y) + y\,V(x,y)] - \frac{2}{\sqrt{\pi}} \tag{22}$$

in conjunction with the relations

$$\frac{\partial L(x,y)}{\partial y} = \frac{\partial V(x,y)}{\partial x} \quad \& \quad \frac{\partial L(x,y)}{\partial x} = -\frac{\partial V(x,y)}{\partial y} \tag{23}$$

Schreier [Schreier 1992] reports numerical problems that can arise when subtracting two numbers of approximately equal magnitudes (when $\partial V/\partial x \sim 0$). Letchworth and Benner [Letchworth et al. 2007] use a special algorithm to calculate these partial derivative (to accuracy <0.5%) but this increases the computational time by about 70%.

## 4- NUMERICAL ANALYSIS AND MACHINE LIMITATIONS

### 4.1. High Accuracy Computations

Due to the finite number of decimal digits available to store a real number in floating-point arithmetic, there are machine limitations on the evaluation of the above summations. The floating-point relative accuracy, $\varepsilon$, and the smallest positive floating-point number, $R_{min}$, in the used computational platform impose restrictions on the accuracy and cause practical machine-truncation of the sums. At any stage during the computation of the sums, the new accumulated sum after adding the term $\alpha_{n+1}$ of the series can be written as

$$\Sigma\big|_{n+1} = \Sigma\big|_n + \alpha_{n+1} = \Sigma\big|_n \times (1 + \Delta(n)), \tag{24}$$

where $\Delta(n) = \alpha_{n+1}/\Sigma\big|_n$.

(24) implies that the sum of $n+1$ terms will not differ from the sum of $n$ terms (i.e. the series will be effectively machine-truncated) if the term $\Delta(n) = \alpha_{n+1}/\Sigma\big|_n$ becomes less than the floating-point relative accuracy, $\varepsilon$, or if $\alpha_{n+1} < R_{min}$. For computational efficiency (shorter computational time) we need to specify these internally truncated terms and exclude them from the computational loop.

Starting with the sum $\Sigma_1$ and considering the possibility of machine-truncation of the sum due to the underflow of the terms $\alpha_{n+1}$, a simple safe estimation for the last value of the index *n* to be included in the evaluation of the sum, $n_{cut}^{\Sigma_1}$, can be derived based on the underflow of the exponential factors only (since the pre-exponential factor is already less than unity)

$$n_{cut,R_{\min}}^{\Sigma_1} \approx \frac{1}{a}\sqrt{-\ln(R_{min})-x^2} \qquad (25)$$

It is implicitly understood that (25) implies rounding to the nearest integer towards infinity. Similarly, the terms of the sums, $\Sigma_2$ and $\Sigma_4$, have the same exponential dependence, while both of the pre-exponential factors will have values less than or equal to unity for $n \geq 1/a$. In such a case, the term $\alpha_{n+1}$ in both sums will reach values $<R_{min}$ if the exponential factor becomes $\leq R_{min}$ which gives

$$n_{cut,R_{\min}}^{\Sigma_2,\Sigma_4} \approx \frac{1}{a}\left[\sqrt{-\ln(R_{min})} - x\right] \qquad (26)$$

Both (25) and (26) indicate that, for $x \geq \sqrt{-\ln(R_{min})}$, no terms from the arguments of the sums $\Sigma_1$, $\Sigma_2$ and $\Sigma_4$ will be effective in the computations and that the values of these sums will be effectively truncated. We note, however, that machine-truncation of these sums due to the $R_{min}$ limitation would also imply machine-truncation due to machine accuracy (i.e, $\Delta(n) \leq \varepsilon$) as seen in (24). For computational efficiency, this later condition may be used to break the computational loop as additional cycles of the computations or additional terms of the series will not change the values of the sums.

Furthermore, we can also make the computations of the sums $\Sigma_3$ and $\Sigma_5$ very efficient. As pointed out above, the values of the terms of these two series grow initially with *n* up to a certain value (peak) and then decay continuously as *n* increases. Simple investigation of these two sums shows that this peak is in the vicinity of *n=x/a*. Accordingly, if one starts calculating these sums from around *n=x/a* and proceeds in both directions, the values of the terms will decrease until they get machine-truncated. For each value of the index *n* used in the calculation of the terms of the sums $\Sigma_1$, $\Sigma_2$ and $\Sigma_4$ we can add to each of the sums $\Sigma_3$ and $\Sigma_5$ two terms by marching one step in each direction. The sums are truncated when the sum of the newly added two terms relative to the value of the previously accumulated sum becomes less than the machine accuracy. Handling the computation of the sums, $\Sigma_3$ and $\Sigma_5$ this way leads to a significant saving in execution time by dramatically reducing the number of terms requiring evaluation which, in turn, leads to a smaller number of loop cycles. The asymptotic expression in the first line of (6) is used for values of $x<R_{min}$.

### 4.2. Accuracy Vs Efficiency Trade-off

With 16-digit floating-point arithmetic and for values of the parameter *a* > ½, the expansion in (8) becomes less accurate and its accuracy will be governed by the corresponding value of the relative error, *E,* as shown in Table 1. It is not necessary, therefore, to keep the strict condition for truncating the sums, i.e., $\Delta(n)<\varepsilon$ since the accuracy of the computations will be governed principally by the relative error *E*. Recalling that the relation between *a* and *E* can be written as $E \sim 2e^{-\pi^2/a^2}$ it seems more appropriate to choose *E* (named *tiny* in the code) instead of *ε* to test for the convergence of the sums. That way we can reduce the number of terms included in the sums by excluding terms that will not effectively enhance the accuracy and thus achieve reasonable acceleration of the computations and reduce the computational time. Noting that, changing the value of *a* changes *E* (*tiny*) and vice versa, we could choose either *a* or *tiny* as the free parameter for controlling the accuracy and efficiency of the computations. This free parameter will be used as an argument of the function. We have chosen to use *tiny* and we calculate *a* internally from the above relation as *tiny* gives a better indication of the accuracy of the computations.

This allows flexibility for accuracy vs efficiency trade-offs while maintaining the ability to run the code for high accuracy.

### 4.3. Calculation of the Exponentials and Other Numerical Considerations

The central part of the present algorithm depends on the evaluation of the sums (15)-(19) which all have exponential terms. The evaluation of intrinsic functions like the exponential function is known to be slower than other simple mathematical operations such as multiplication and/or division. A naïve computation of these sums would require three exponential evaluations per computational loop. This would be computationally expensive. However, we may reduce this to just one exponential evaluation in each cycle of the loop. For the case of $x < \sqrt{-\ln(R_{min})}$ we have some flexibility in evaluating the term $e^{-x^2}$ either separately or by combining it with other terms. In such a case one can write all of the exponentials in (15)-(19) in terms of the exponential $e^{-a^2n^2}$ where

$$e^{-(a^2n^2+x^2)} = e^{-a^2n^2} \times e^{-x^2} \qquad (27)$$

$$e^{-(an+x)^2} = e^{-a^2n^2} \times e^{-x^2} \times \prod_1^n e^{-2ax} \qquad (28)$$

$$e^{-(an-x)^2} = e^{-a^2n^2} \times e^{-x^2} \times \prod_1^n e^{2ax} \qquad (29)$$

In the above expressions $e^{-x^2}$, $e^{-2ax}$ and $e^{2ax}$ are calculated once outside the loop, for each value of *x,* and the products are simply performed using multiplies inside the computational loop. For $x \geq \sqrt{-\ln(R_{min})}$ only $\Sigma_3$ and $\Sigma_5$ (which have the same exponential factor) contribute to the calculation of the real and imaginary parts of the Faddeyeva function. However, due to the nature of the computation of these terms, the exponential factor needs to be computed twice; once for the step to the right of $n_0$ =ceil(*x/a*) and once on the left wing where "ceil" indicates rounding to the nearest integer towards infinity. The indices for these two factors are related to the loop index, *n*, by $n_{3\text{-}plus} = n_0+(n\text{-}1)$ and $n_{3\text{-}minus} = n_0 - n$ if $n_{3\text{-}minus} \geq 1$, respectively. When $n_{3\text{-}minus} <1$, only the term on the right wing is included. The exponential factors for the terms on the left and right wings are related by

$$e^{-(an_{3\text{-}minus}-x)^2} = e^{-(an_{3\text{-}plus}-x)^2} \times e^{(a^2+2ax-2a^2n_0)} \times \prod_1^n e^{(4a^2n_0-4ax-2a^2)} \qquad (30)$$

Clearly the second exponential factor on the right hand side of (30) and the argument of the product may be calculated once outside the loop, for each value of *x*, thus reducing the number of exponential function evaluations to only one per cycle. The product can be evaluated using just multiplies inside the loop, thus reducing the computational time.

A few more important computational points are related to the calculation of the imaginary part of the Faddeyeva function using (14). Firstly, for values of $x<<1$, the two sums $\Sigma_4$ and $\Sigma_5$ become very close to each other and the subtraction of these two sums could significantly affect the accuracy in regions of the computational domain where these two terms are the dominant terms in calculating the imaginary part of the Faddeyeva function. However, this problem can be simply overcome by expressing the sum of these two terms (for $x<<1$) in its original form as in (12), i.e. in terms of $\sinh(2anx)$. The first three terms in the series expansion of $\sinh(x)$ will be sufficient to express $\sinh(x)$ to the machine accuracy for $x \leq 10^{-2}$ and since $2an$ is usually $<20$ in the present computations then this is satisfactorily for $x \leq 5 \times 10^{-4}$.

The second important point in the calculation of the imaginary part of the Faddeyeva function using (14) is related to the calculation of the sum of the first three terms on the right hand side of the equation, for $x < \sqrt{-\ln(R_{min})}$, which can be written safely in the form

$$e^{-x^2}\sin(2xy)\left\{-\operatorname{erfcx}(y) + \frac{a}{y\pi}\left[1 + 2\sum_{n=1}^{\infty}\frac{e^{-a^2n^2}}{1+a^2n^2/y^2}\right]\right\} \qquad (31)$$

The terms in the curly brackets are only dependent on y and for $y \geq 5$ we have found that this sum is zero to machine accuracy. Using this prevents rounding errors affecting the accuracy of computations. Note that for very small values of $x$ the result of the whole expression (31) is O($yx$) while the total of $-\Sigma_4 + \Sigma_5$ is O($x$), the significance of rounding errors is thus clear for small values of $x$ and relatively large values of $y$.

## 5. THE MATLAB FUNCTION FADDEYEVA.M

The function Faddeyeva(**z**,*tiny*) returns, in general, an array of complex values for the Faddeyeva function of the same size as the input array for the complex variable $z$. The input **z** is usually an array (with one or two dimensions) but can be a single scalar as well. When **z** contains only imaginary values **z**=i**y,** the function returns the real values calculated from the MATLAB built-in function "erfcx(**y**)". The function is set for the calculation for the whole complex domain. However, for negative values of **y** and exp(**y**$^2$-**x**$^2$) greater than the largest floating point number in the computational platform, Faddeyeva cannot calculate the Faddeyeva function due to inescapable overflow problems. The function checks for acceptable values and issues an error message for any points outside this domain.

The value of the scalar free parameter "*tiny"* can be chosen by the user within the range $tiny_{min} \leq tiny \leq 10^{-4}$ to control the accuracy and computational time. The value of $tiny_{min}$ is a value close to but less than the floating-point relative accuracy, $\varepsilon$. For example, for a 16 digit computational platform, $tiny_{min}$ can be taken roughly to be $\sim 1.43 \times 10^{-17}$ (the value of $E$ corresponding to $a=1/2$ in Table 1) while for a 32-digit computational platform $tiny_{min}$ can be taken to be roughly $10^{-33}$. The maximum value of $tiny=10^{-4}$ corresponds to $a=1$ (the maximum value for $a$ for which the expansion in (8) can be used). Increasing the value of *tiny* within its above mentioned range will decrease the computational time at the expense of the computational accuracy and vice versa.

Choosing a value of $tiny<tiny_{min}$ will just increase the run time without any improvement in the accuracy of computations which will then be governed solely by the machine characteristics. Values of *tiny*<*tiny*$_{min}$ or *tiny*>$10^{-4}$ result in *tiny* being reset to *tiny*$_{min}$ or $10^{-4}$ respectively; a warning message is returned in both cases.

It is to be noted that *tiny*<*ε* is used only for the calculation of the corresponding value of the parameter *a* and not for the truncation of the sums since the calculations cannot be claimed to be performed for relative accuracy less than the machine accuracy epsilon, *ε*, in any case. Accordingly, for the truncation of the sums, the maximum of *tiny* and *ε* is used.

## 6. ALGORITHM VERIFICATION AND EFFICIENCY

### 6.1. High Accuracy Computations

Three different independent computational techniques are used to investigate the accuracy of the present algorithm

1- Mathematica™ [Wolfram Research, Inc. 2008] provides the imaginary error function erfi(*z*) as a special function which can be evaluated and then used in conjunction with the relation given in (1); that is $w(z) = e^{-z^2}(1 + i\,\text{erfi}(z))$, to calculate the Faddeyeva function [Zaghloul 2008]. The arbitrary-precision arithmetic used in Mathematica™ allows us to obtain highly accurate values for the function erfi(*z*) although these calculations are very expensive computationally. This is only generally suitable for applications where the speed of arithmetic is not a restrictive factor, or where precise results for a small number of evaluations are required. We can, however, generate highly accurate values of the function erfi(*z*) using Mathematica™ by using large numbers of digits of precision,

2- the simple proper integral given in reference [Zaghloul 2007] can be used to calculate the real Voigt function (real part of the Faddeyeva function), and

3- the Algorithm 680 [Poppe et al. 1990], which is widely used in the literature and is implemented in many software packages and libraries, calculates the Faddeyeva function to a claimed accuracy of 14 significant digits.

The relatively long computational time associated with the first two methods makes them inefficient for use in applications requiring a large number of function evaluations. For this reason, Algorithm 680 is regarded as the competitive highly accurate algorithm due to its computational speed and claimed accuracy.

Table 2 presents sample results calculated using these three methods with the results from the present algorithm. The values of the complex variable *z* used in the computations in this table have been selected to allow some conclusions to be drawn in addition to establishing confidence and reliability in the present code. The value of the relative error in the calculation of the Faddeyeva function proposed here is given in Table 3 and compared with the other approaches, taking the values of $w$ calculated using the function erfi(*z*) from Mathematica™ with high number of digits of accuracy as reference values.

Looking closer into the values given in these two tables we conclude that;
a) compared to calculating the Faddeyeva function using erfi(*z*) from Mathematica™, the present algorithm is more reliable since we failed to calculate erfi(*z*) using Mathematica™ for values of *x*>$3.9\times10^4$ and *y*>$2.8\times10^4$ while the present algorithm does not suffer such a limitation. Note that for this domain Mathematica™ cannot be used as reference, and therefore, no comparison was performed for this range in Table 3.

b) for the whole domain of computations, the present algorithm shows very high accuracy as shown in Table 3, while the Algorithm 680 suffers a catastrophic loss of accuracy in the vicinity of *x*=6.3 as well as for small values of *y* signified by *bold-face* numbers in Tables 2 and 3. The relative error for the real part of the Faddeyeva function from Algorithm 680 in this region of the first quadrant goes up to 100%. It has to be noted that *x*=6.3 is one of the built-in values in Algorithm 680.

To investigate the effectiveness of the present algorithm compared to Algorithm 680, we calculated the Faddeyeva function using both algorithms for 2,840,710 points of the complex variable *z* distributed over the upper half of the complex plane using the grid *y=logspace(-20, 4, 71)* and *x=linspace(-200, 200, 40001)* where *y=logspace(-20, 4, 71)* generates a row vector of 71 logarithmically equally spaced points between $10^{-20}$ and $10^4$ and *x=linspace(-200, 200, 40001)* generates a row vector of 40001 linearly equally spaced points between -200 and 200. Using *Matlab 7.9.0.529 (R2009b)*, the computational time taken by the present algorithm was found to be <8.0% of that taken by Algorithm 680, which represents a significant time saving.

Figures (1-a) and (1-b) show surface plots of the absolute relative error $\delta_V = |V - V_{ref}|/V_{ref}$ and $\delta_L = |L - L_{ref}|/L_{ref}$ in the results obtained from the present algorithm using results from Algorithm 680, which is available in Matlab, as reference values. As can be clearly seen from these figures, the results from the present algorithm show high agreement (around 13 significant digits) over the chosen computational domain except for the region in the vicinity of *x*=6.3 and small values of *y* where Algorithm 680 badly loses its accuracy as indicated above. The above findings provide the necessary verification and confirm the high accuracy as well as the reliability of the present algorithm.

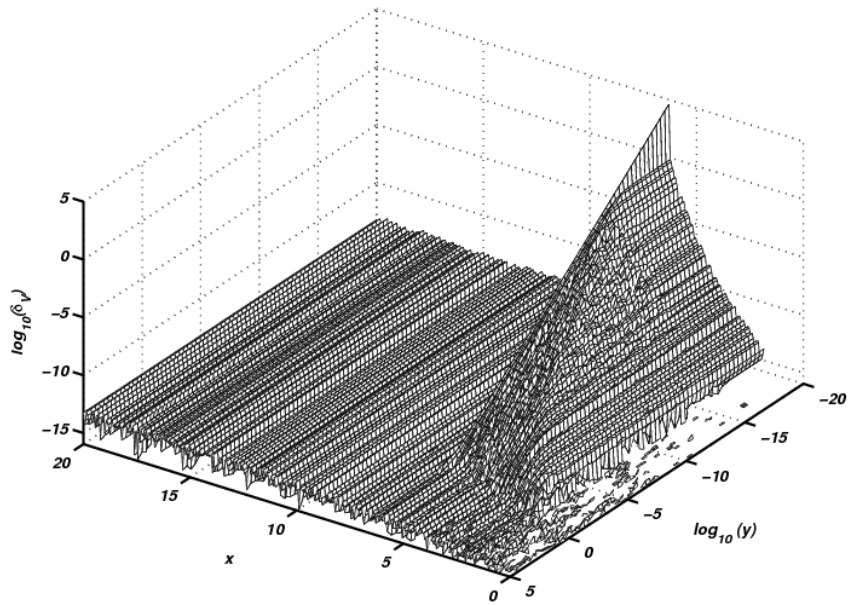

**Figure 1-a** Absolute relative error $\delta_V = |V - V_{ref}|/V_{ref}$ in the calculations of the real part of the Faddeyeva function from the present algorithm using the results from Algorithm 680 as reference values.

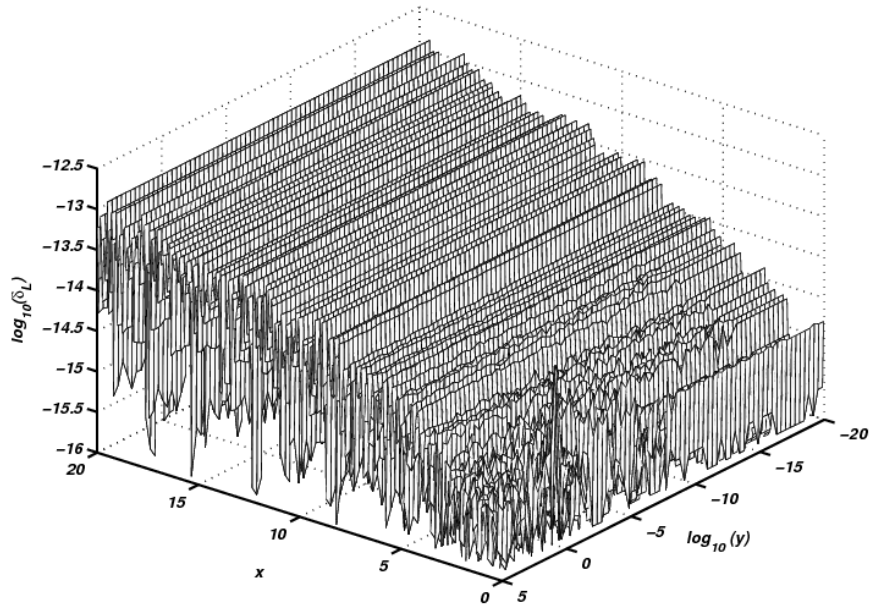

**Figure 1-b** Absolute relative error $\delta_L = |L - L_{ref}|/L_{ref}$ in the calculations of the imaginary part of the Faddeyeva function from the present algorithm using the results from Algorithm 680 as reference values.

## 6.2. Efficient Computations with Lower Accuracies

Table 4 below shows values of the free variable *tiny* used in the calling argument of our Matlab function and the corresponding relative accuracy $\delta_V = |(V - V_{ref})/V_{ref}|$ and $\delta_L = |(L - L_{ref})/L_{ref}|$ in the calculations, using values calculated with the highest accuracy obtainable from the present algorithm as reference values. The need to quantify the efficiency improvements obtainable when using the accuracy vs efficiency trade-off capability of the present algorithm is the reason of using the highest accuracy computations from the present algorithm as reference values. The run times required to calculate the function for 2,840,710 points generated using the grid *y=logspace*(-20, 4, 71) and *x=linspace*(-200, 200, 40001) relative to the run time required to perform the same computations using the highest accuracy computations from the present algorithm are also included in the table. As can be seen from the table, running the present algorithm at lower accuracy improves the efficiency of the computations and decreases the computational time by up to 45%. Compared to other efficient and low-accuracy algorithms in the literature [Hui et al. 1978; Humlíček 1982; Poppe et al. 1990; Lether et al. 1991; Shippony et al. 1993; Weideman 1994], the present algorithm seems to be more reliable even at low-accuracy. In addition, other algorithms fail in some regions of the computational domain, particularly near the real axis (very small values of *y*); for example, the Poppe and Wijers algorithm [Poppe et al. 1990], known for its accuracy, fails in this region, returning results for the real part of the Faddeyeva function that are several orders of magnitude away from the correct values.

Figure 2 shows a comparison between the calculations of the partial derivative $\partial V(x,y)/\partial x$ using the present algorithm (run at the lowest accuracy) and calculations from Algorithm 680, for $y=10^{-20}$, in the region *x*=[6.1-6.5]. Here we see that the results from the present algorithm (even when run at the lowest accuracy) seem to be more accurate and more reliable than computations from Poppe and Wijers algorithm in this region, where the latter loses its accuracy and fails to produce the correct behavior of $\partial V(x,y)/\partial x$.

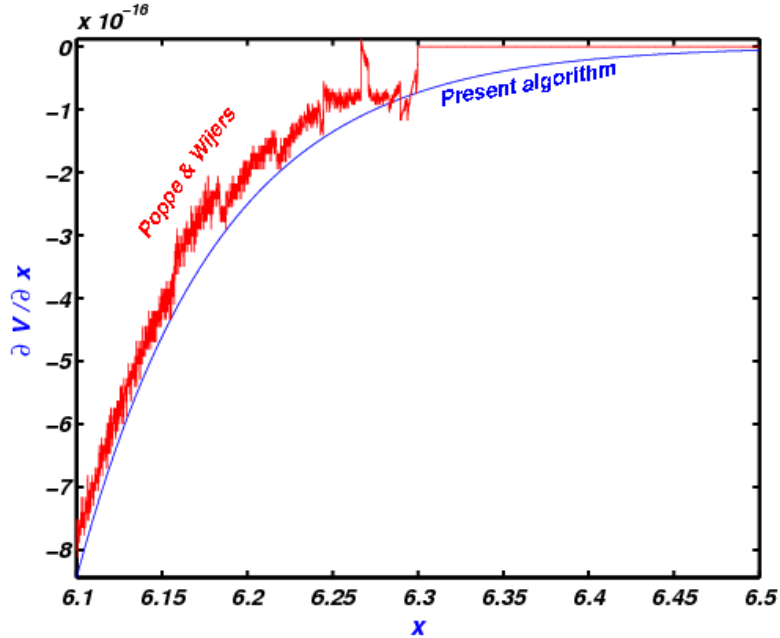

**Figure 2** $\partial V(x,y)/\partial x$ **as calculated, from the present algorithm (*tiny*=10$^{-4}$) and from Poppe and Wijers algorithm, using (21), for *y*=10$^{-20}$.**

Figure 3 shows a surface plot of the relative error $\delta_V = \left|(V - V_{ref})/V_{ref}\right|$ for the results obtained from Hui's algorithm [Hui et al. 1978] using results from the present algorithm as reference values. We note that the Matlab version of Hui's algorithm (cerf.m [Hui et al. 1978]) employed in this comparison, uses the *p*=5 rational approximation where *p* is the degree of nominator polynomial. As is clear from the plot, the relative errors in the results of Hui's algorithm are very large for small values of *y* and reach 14 order of magnitude for medium values of *x* when *y*=10$^{-20}$.

In addition to the large errors for medium values of *x* and small values of *y*, Hui's algorithm produces negative values for the Voigt function (real part of the Faddeyeva function) for example, for *y*=10$^{-5}$ and *x*=4. The Voigt function is positive over the whole first quadrant.

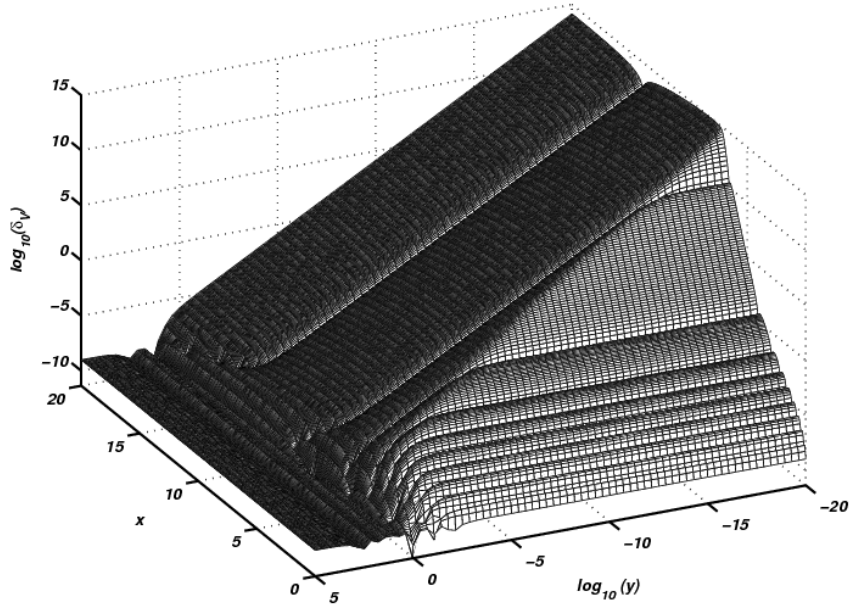

**Figure 3** Absolute relative error, $\delta_V = \left|(V - V_{ref})/V_{ref}\right|$, in the calculations of the real part of the Faddeyeva function from Hui's algorithm using results from the present algorithm (*tiny=tiny$_{min}$*) as reference values.

Figure 4 shows a comparison between the calculations of the partial derivative $\partial V(x,y)/\partial x$ using the present algorithm (run at the lowest accuracy) and calculations from Hui's algorithm, using (21), for $y=10^{-20}$, in the region $x=[7,15]$. As we see from the figure, calculations from the present algorithm seem to be more accurate and more reliable than Hui's algorithm which fails to produce the correct behavior (negative $\partial V(x,y)/\partial x$ ) or the correct order of magnitude of $\partial V(x,y)/\partial x$ in this region of the computational domain.

We note that the error contours given in [Hui et al. 1978] were presented either for the *modulus* of the complex error function or for the *absolute* value of the Voigt function $V(x,y)$. That was, probably, the reason why such failures were not clear from their paper.

Although Hui's algorithm takes about 10% of the computational time taken by the present algorithm (for *tiny*=$10^{-4}$) and about 5% of the computational time taken by the present algorithm for the highest-accuracy computations, the fact that it fails to produce the correct values or correct signs of the function or even its correct behavior, in this region of the computational domain, poses important questions about its reliability.

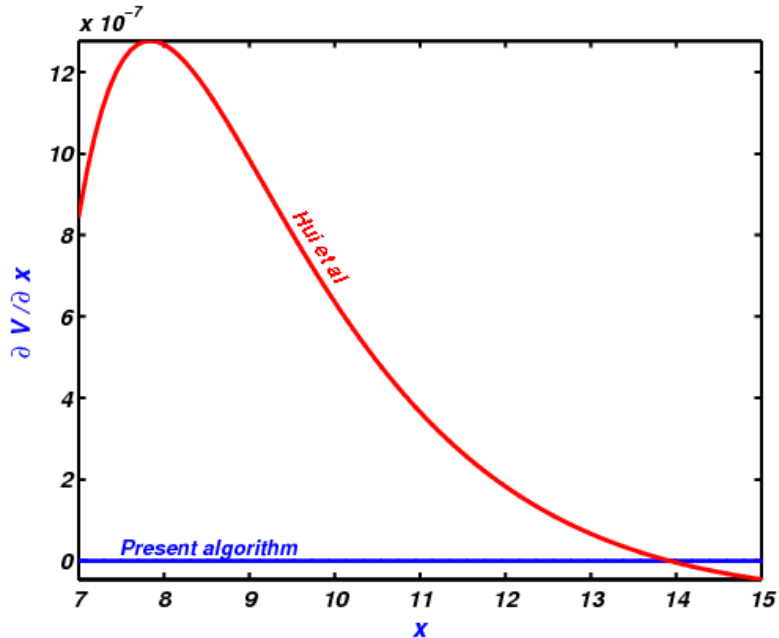

**Figure 4** $\partial V(x,y)/\partial x$ **as calculated, from the present algorithm (*tiny*=10$^{-4}$) and from Hui's algorithm, using (21), for *y*=10$^{-20}$.**

Humlíček [Humlíček 1982] reported that any rational approximation suffers inevitable failure near the real axis and he attempted to overcome this failure in his algorithm, w4. However, investigating the results of Humlíček's algorithm we found that it also suffers complete failure near the real axis in the vicinity of *x*=5.5. The results for *y*=10$^{-20}$ show that Humlíček's algorithm underestimates the real part of the Faddeyeva function by 8 orders of magnitude and by 3 orders of magnitude for *y*=10$^{-15}$. Figure 5 shows a surface plot of the relative error in the calculation of *V(x,y)* using Humlíček's algorithm taking the results from the present algorithm as reference.

Table 5 shows that the computational time using the Humlíček's original code is almost three times that used by the present algorithm (for the highest-accuracy computations). Even using a more efficient version of Humlíček's code, modified by the present authors, we note that the computational time taken by Humlíček's algorithm is still longer than that taken by the present algorithm (for the highest-accuracy computations).

Figure 6, on the other hand, shows a comparison between the calculations of the partial derivative $\partial V(x,y)/\partial x$ using (21) from the present algorithm (run at the lowest accuracy) and those calculations from Humlíček's algorithm, for *y*=10$^{-20}$, in the region *x*=[5.4,6.4]. This shows that Humlíček's algorithm does not produce the correct behavior of $\partial V(x,y)/\partial x$ in this domain.

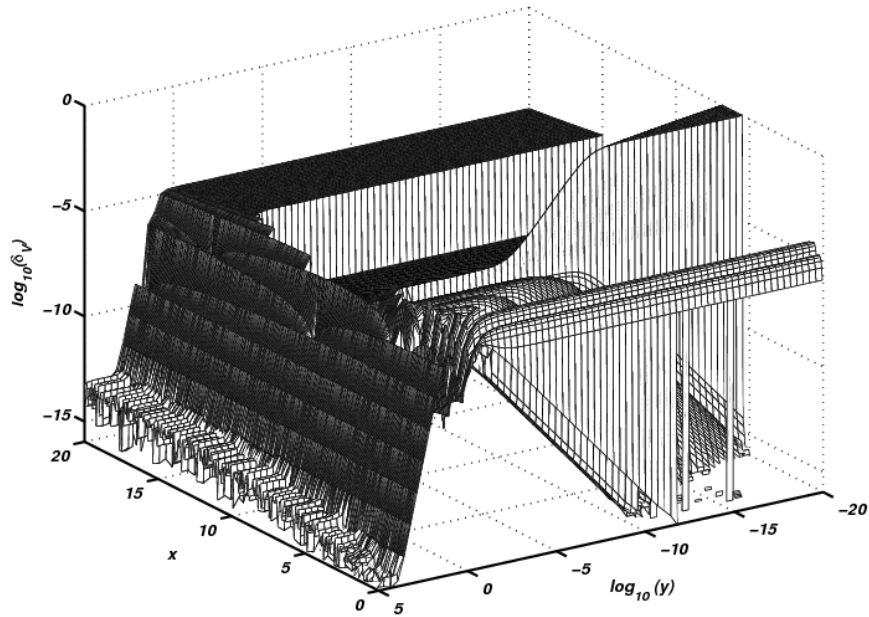

**Figure 5** Absolute relative error, $\delta_V = \left|(V - V_{ref})/V_{ref}\right|$, in the calculations of the real part of the Faddeyeva function from Humlíček's algorithm using results from the present algorithm (*tiny=tiny$_{min}$*) as reference values.

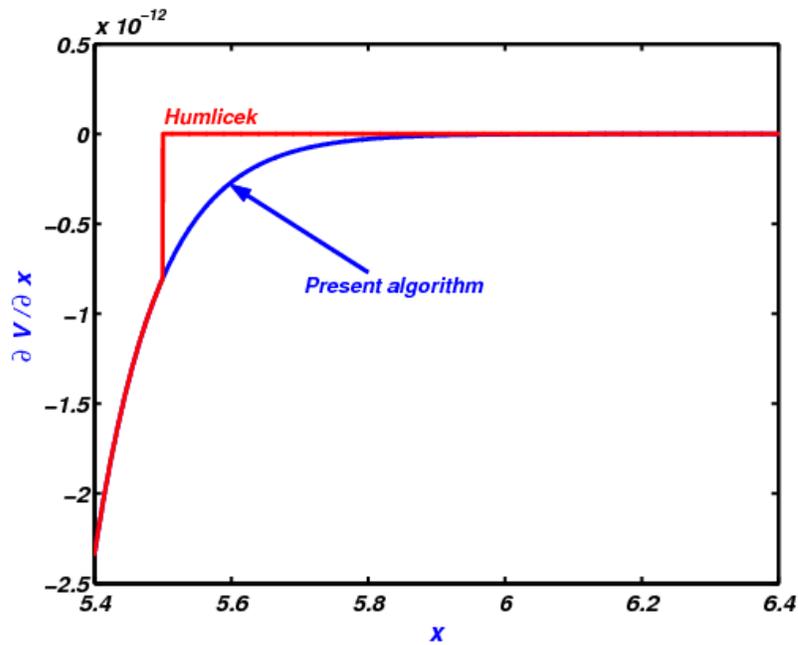

**Figure 6** $\partial V(x,y)/\partial x$ as calculated, from the present algorithm (*tiny*=$10^{-4}$) and from Humlíček's algorithm, using (21), for $y=10^{-20}$.

As with Hui's algorithm, Weideman's algorithm [Weideman 1994] also produces negative values for the real part of the Faddeyeva function near the real axis. The negative values for the Voigt function calculated from Weideman's algorithm appear for all values of the parameter N (number of terms in the rational series) for $y=10^{-20}$. A surface plot of the relative error in the calculations of the real part of the Faddeyeva function using Weideman's algorithm with N=256 taking the results of the present algorithm as a reference is shown in Figure 7. The figure shows that the errors resulting from Weideman's algorithm are catastrophic for small values of $y$ and that the computed magnitude of $V(x,y)$ is overestimated by up to 6 orders of magnitude for $y=10^{-20}$. The situation becomes even worse for smaller values of N. Table 5 shows that the run time of the present algorithm at highest accuracy is shorter than the run time of Weideman's algorithm with N=256.

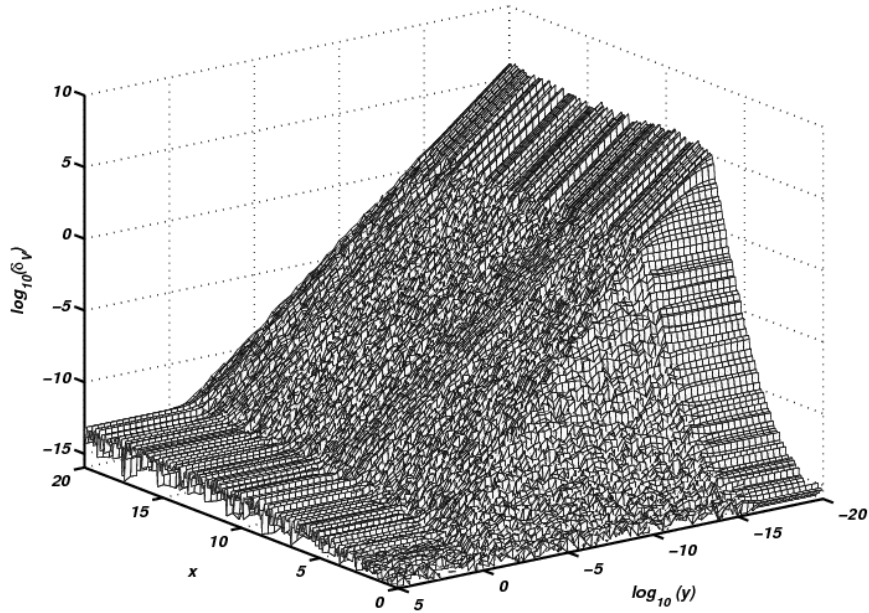

**Figure 7 Absolute relative error, $\delta_V = \left|(V - V_{ref})/V_{ref}\right|$, in the calculations of the real part of the Faddeyeva function from Weideman's algorithm, N=256, using results from the present algorithm (*tiny=tiny$_{min}$*) as reference values.**

Figures (8-a) and (8-b) show comparisons between the calculations of the partial derivative $\partial V(x,y)/\partial x$ using the present algorithm (run at the lowest accuracy) and Weideman's algorithm with N=128 and N=32, respectively. The calculations shown in the figure are for the region $x=[6,15]$ and $y=10^{-20}$ in Figure (8-a) and $y=10^{-10}$ in Figure (8-b). From these figures, we see that Weideman's algorithm does not reproduce the correct behavior for $\partial V(x,y)/\partial x$ in the regions shown.

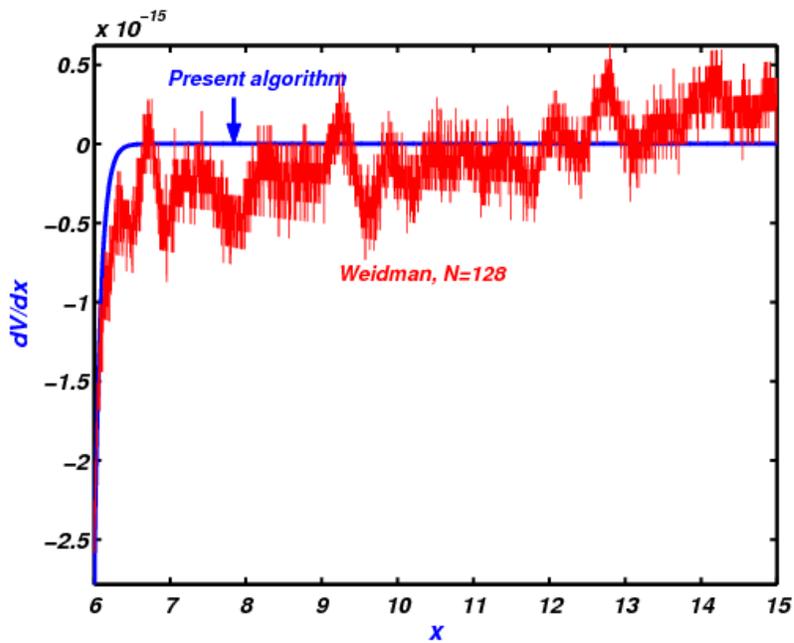

**Figure 8-a** $\partial V(x,y)/\partial x$ as calculated, from the present algorithm (*tiny*=10⁻⁴) and from Weideman's algorithm with N=128, using (21), for $y=10^{-20}$.

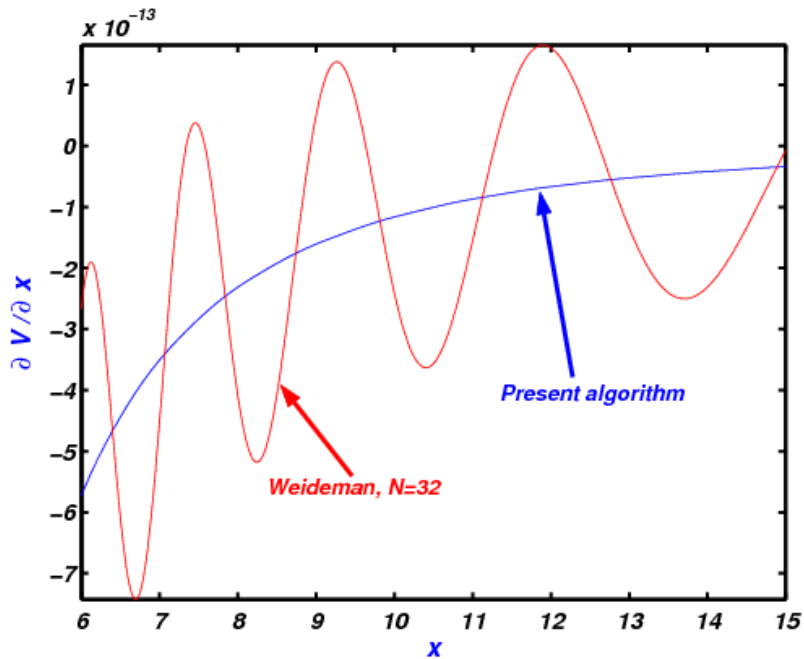

**Figure 8-b** $\partial V(x,y)/\partial x$ as calculated, from the present algorithm (*tiny*=10⁻⁴) and from Weideman's algorithm with N=32, using (21), for $y=10^{-20}$.

Other competitive algorithms in the literature also show loss of accuracy in some regions of the computational domain. For example, the algorithm by Shippony and Read [Shippony et al. 1993] exhibits the same failure in calculating the real part of the Faddeyeva function near the real axis. In particular we detected the same loss of accuracy suffered by Poppe and Wijers algorithm for very small values of $y$ near $x$=6.3. In addition, the Shippony and Read algorithm produces negative values for $V(x,y)$ and/or $L(x,y)$ in several regions of the computational domain. Just for example we refer to the points at $x$=1.5 and $y$=1.5, 2.0, 2.5, 3.0, 3.5, 4.0, 5.0 etc. It has to be noted that these failures have been obtained even with the use of the correction provided by Shippony and Read in [Shippony et al. 2003].

The situation is no better with the Letchworth and Benner's algorithm [Letchworth et al. 2007] where similar failures and loss of accuracy are obtained for very small values of $y$ and values of $x$ greater than but close to $x$=5.76. For $x$=5.76 and $y=10^{-20}$ Letchwoth and Benner's algorithm returns, for the real part of Faddeyeva function the value $V$=1.783900323491466×$10^{-22}$ while the value returned from the present algorithm is $V$=3.900779639194698×$10^{-15}$ and that returned from using the function erfi($z$) from Mathematica$^{TM}$ is 3.900779639194697×$10^{-15}$. The algorithm in [Zaghloul 2007] returns 3.900779639194698×$10^{-15}$. Figure 9 shows a similar comparison between the calculations of the partial derivative $\partial V(x,y)/\partial x$ using the present algorithm (with $tiny$=$10^{-4}$) and calculations from Letchworth and Benner algorithm, in (21), for $y=10^{-20}$. The failure of Letchworth and Benner's algorithm for values of $x$ greater than but close to $x$=5.76 can be easily recognized from the figure.

It has to be emphasized here that these "competitive codes" were (probably) not designed for such extreme values of $y$ and the tests presented here are more a demonstration of the "superior accuracy" of our software even for extreme values.

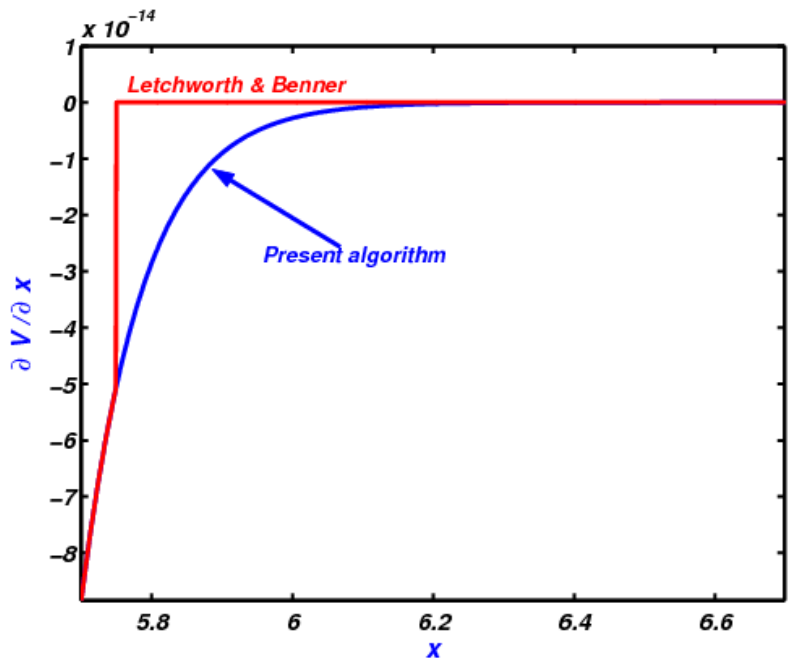

**Figure 9** $\partial V(x,y)/\partial x$ as calculated, from the present algorithm ($tiny$=$10^{-4}$) and from Letchworth & Benner's algorithm, using (21), for $y=10^{-20}$.

## 7. CONCLUSIONS

An algorithm accompanied by a computer code, in the form of a MATLAB$^{TM}$ function, for the numerical evaluation of the Faddeyeva function $w(z)$ is presented. The algorithm is more accurate and avoids failures discovered in other competitive published algorithms. In addition to its superior accuracy, the present algorithm and computer code allow a flexible accuracy vs efficiency trade-off through controlling a free parameter *tiny*. By adjusting the value of this parameter the function can be run with a lower accuracy and shorter computational time or high accuracy and longer computational time. Even when run at its lowest accuracy, the present algorithm avoids major problems suffered by other competitive codes. For all levels of accuracy the present code is safer and more reliable since it does not return negative values for real and/or imaginary parts of the Faddeyeva function nor does it suffer from the loss of accuracy exhibited by some of the other competitive codes. The present algorithm can, therefore, be safely used and implemented in personal and commercial libraries.

Table 2: Results from algorithms in the literature in comparison with results from the present algorithm for some selected values of $z$

| $z$ | | Mathematica$^{TM}$ | | Algorithm 680 | | [Zaghloul 2007] | Present algorithm | |
|---|---|---|---|---|---|---|---|---|
| $x$ | $y$ | $V$ | $L$ | $V$ | $L$ | $V$ | $V$ | $L$ |
| 6.3e-002 | 1.0e-020 | 9.960388660702479e-001 | 7.090008726353683e-002 | 9.960388660702479e-001 | 7.090008726353685e-002 | 9.960388660702479e-001 | 9.960388660702479e-001 | 7.090008726353669e-002 |
| 6.3e-002 | 1.0e-014 | 9.960388660702367e-001 | 7.090008726353558e-002 | 9.960388660702367e-001 | 7.090008726353560e-002 | 9.960388660702366e-001 | 9.960388660702366e-001 | 7.090008726353543e-002 |
| 6.3e-002 | 1.0e-012 | 9.960388660691284e-001 | 7.090008726341133e-002 | 9.960388660691284e-001 | 7.090008726341135e-002 | 9.960388660691284e-001 | 9.960388660691284e-001 | 7.090008726341118e-002 |
| 6.3e-002 | 1.0e-010 | 9.960388659583033e-001 | 7.090008725098674e-002 | 9.960388659583034e-001 | 7.090008725098676e-002 | 9.960388659583034e-001 | 9.960388659583033e-001 | 7.090008725098659e-002 |
| 6.3e-002 | 1.0e-006 | 9.960377466254799e-001 | 7.089996176278113e-002 | 9.960377466254800e-001 | 7.089996176278113e-002 | 9.960377466254802e-001 | 9.960377466254801e-001 | 7.089996176278086e-002 |
| 6.3e-002 | 1.0e-002 | 9.849424862549036e-001 | 6.965909657459020e-002 | 9.849424862549037e-001 | 6.965909657459021e-002 | 9.849424862549038e-001 | 9.849424862549039e-001 | 6.965909657459005e-002 |
| 6.3e-002 | 1.0e+001 | 5.613881832823887e-002 | 3.502232333332985e-004 | 5.613881832823888e-002 | 3.502232333332986e-004 | 5.613881832823887e-002 | 5.613881832823886e-002 | 3.502232333332973e-004 |
| 6.3e-002 | 1.2e+001 | 4.685295149211636e-002 | 2.442987772965768e-004 | 4.685295149211637e-002 | 2.442987772965769e-004 | 4.685295149211637e-002 | 4.685295149211637e-002 | 2.442987772965766e-004 |
| 6.3e-002 | 1.5e+001 | 3.752895161491573e-002 | 1.569287266610685e-004 | 3.752895161491574e-002 | 1.569287266610686e-004 | 3.752895161491573e-002 | 3.752895161491574e-002 | 1.569287266610681e-004 |
| 6.3e-002 | 2.0e+002 | 2.820912377324508e-003 | 8.885651855627418e-007 | 2.820912377324509e-003 | 8.885651855627419e-007 | 2.820912377324511e-003 | 2.820912377324508e-003 | 8.885651855627396e-007 |
| 6.3e-002 | 1.0e+005 | *Fails to evaluate erfi* | | 5.641895835193230e-006 | 3.554394375816296e-012 | 5.641895835475324e-006 | 5.641895835193228e-006 | 3.554394375816285e-012 |
| | | | | | | | | |
| 6.3e+000 | 1.0e-020 | **5.792460778844102e-018** | 9.072765968412736e-002 | **1.478934492449188e-022** | 9.072765968412733e-002 | **5.792312885394871e-018**$^{\$}$ | **5.792460778844116e-018** | 9.072765968412679e-002 |
| 6.3e+000 | 1.0e-014 | **1.536857621303171e-016** | 9.072765968412736e-002 | **1.478934492449189e-016** | 9.072765968412733e-002 | **1.536857621303171e-016** | **1.536857621303163e-016** | 9.072765968412679e-002 |
| 6.3e+000 | 1.0e-012 | **1.479513723737762e-014** | 9.072765968412736e-002 | **1.478934492449189e-014** | 9.072765968412733e-002 | **1.479513723737762e-014** | **1.479513723737753e-014** | 9.072765968412679e-002 |
| 6.3e+000 | 1.0e-010 | **1.478940284762108e-012** | 9.072765968412736e-002 | **1.478934492449189e-012** | 9.072765968412733e-002 | **1.478940284762108e-012** | **1.478940284762099e-012** | 9.072765968412679e-002 |
| 6.3e+000 | 1.0e-006 | **1.478934493028413e-008** | 9.072765968412492e-002 | **1.478934492449148e-008** | 9.072765968412488e-002 | **1.478934493028413e-008** | **1.478934493028404e-008** | 9.072765968412433e-002 |
| 6.3e+000 | 1.0e-002 | 1.478930389133942e-004 | 9.072741516349275e-002 | 1.478930389133851e-004 | 9.072741516349273e-002 | 1.478930389133943e-004 | 1.478930389133934e-004 | 9.072741516349218e-002 |
| 6.3e+000 | 1.0e+001 | *4.040671157393860e-002* | *2.527577277549421e-002* | 4.040671157393860e-002 | 2.527577277549422e-002 | 4.040671157393859e-002 | 4.040671157393835e-002 | 2.527577277549405e-002 |
| 6.3e+000 | 1.2e+001 | *3.684277239564821e-002* | *1.923808857910893e-002* | 3.684277239564821e-002 | 1.923808857910892e-002 | 3.684277239564818e-002 | 3.684277239564798e-002 | 1.923808857910881e-002 |
| 6.3e+000 | 1.5e+001 | *3.194834330452624e-002* | *1.336797114261604e-002* | 3.194834330452625e-002 | 1.336797114261605e-002 | 3.194834330452623e-002 | 3.194834330452605e-002 | 1.336797114261596e-002 |
| 6.3e+000 | 2.0e+002 | *2.818116555672224e-003* | *8.876845457496914e-005* | 2.818116555672223e-003 | 8.876845457496911e-005 | 2.818116225691620e-003 | 2.818116555672206e-003 | 8.876845457496856e-005 |
| 6.3e+000 | 1.0e+005 | *Fails to evaluate erfi* | | 5.641895812802784e-006 | 3.554394361710315e-010 | 5.641895813084879e-006 | 5.641895812802746e-006 | 3.554394361710292e-010 |
| | | | | | | | | |
| 6.3e+002 | 1.0e-020 | 1.421495882582394e-026 | 8.955401496757104e-004 | 1.421495882582395e-026 | 8.955401496757104e-004 | 1.421490510324405e-026 | 1.421495882582395e-026 | 8.955401496757105e-004 |
| 6.3e+002 | 1.0e-014 | 1.421495882582394e-020 | 8.955401496757104e-004 | 1.421495882582395e-020 | 8.955401496757104e-004 | 1.421490510324405e-020 | 1.421495882582395e-020 | 8.955401496757105e-004 |
| 6.3e+002 | 1.0e-012 | 1.421495882582394e-018 | 8.955401496757104e-004 | 1.421495882582395e-018 | 8.955401496757104e-004 | 1.421490510324405e-018 | 1.421495882582395e-018 | 8.955401496757105e-004 |
| 6.3e+002 | 1.0e-010 | 1.421495882582394e-016 | 8.955401496757104e-004 | 1.421495882582395e-016 | 8.955401496757104e-004 | 1.421490510324405e-016 | 1.421495882582395e-016 | 8.955401496757105e-004 |
| 6.3e+002 | 1.0e-006 | 1.421495882582394e-012 | 8.955401496757104e-004 | 1.421495882582394e-012 | 8.955401496757104e-004 | 1.421490510324405e-012 | 1.421495882582395e-012 | 8.955401496757105e-004 |
| 6.3e+002 | 1.0e-002 | 1.421495882224241e-008 | 8.955401494500753e-004 | 1.421495882224242e-008 | 8.955401494500753e-004 | 1.421490509966257e-008 | 1.421495882224242e-008 | 8.955401494500755e-004 |
| 6.3e+002 | 1.0e+001 | 1.421137820009847e-005 | 8.953145713915760e-004 | 1.421137820009848e-005 | 8.953145713915762e-004 | 1.421132452261351e-005 | 1.421137820009847e-005 | 8.953145713915762e-004 |
| 6.3e+002 | 1.2e+001 | 1.705176395541706e-005 | 8.952153529445874e-004 | 1.705176395541707e-005 | 8.952153529445876e-004 | 1.705169956622711e-005 | 1.705176395541707e-005 | 8.952153529445874e-004 |
| 6.3e+002 | 1.5e+001 | 2.131035743074597e-005 | 8.950327582962093e-004 | 2.131035743074598e-005 | 8.950327582962094e-004 | 2.131027699897097e-005 | 2.131035743074598e-005 | 8.950327582962094e-004 |
| 6.3e+002 | 2.0e+002 | 2.582702147491469e-004 | 8.135493143556982e-004 | 2.582702147491469e-004 | 8.135493143556982e-004 | 2.582694362772975e-004 | 2.582702147491469e-004 | 8.135493143556982e-004 |
| 6.3e+002 | 1.0e+005 | *Fails to evaluate erfi* | | 5.641671917237129e-006 | 3.554253307503979e-008 | 5.641671917519157e-006 | 5.641671917237128e-006 | 3.554253307503980e-008 |
| | | | | | | | | |
| 1.0e+000 | 1.0e-020 | 3.678794411714423e-001 | 6.071577058413937e-001 | 3.678794411714423e-001 | 6.071577058413937e-001 | 3.678794411714423e-001 | 3.678794411714423e-001 | 6.071577058413938e-001 |
| 5.5e+000 | 1.0e-014 | 7.307386729528773e-014 | 1.043674364367812e-001 | 7.308245082486227e-014 | 1.043674364367812e-001 | 7.307386729528773e-014$^{*}$ | 7.307386729528773e-014 | 1.043674364367812e-001 |
| 3.9e+004 | 1.0e+000 | *Fails to evaluate erfi* | | 3.709333226385424e-010 | 1.446639957339204e-005 | 3.709333222727304e-010 | 3.709333226385423e-010 | 1.446639957339204e-005 |
| 1.0e+000 | 2.8e+004 | *Fails to evaluate erfi* | | 2.014962794529686e-005 | 7.196295685569928e-010 | 2.014962795814739e-005 | 2.014962794529686e-005 | 7.196295685569929e-010 |

* This is the correct value as calculated using (4) in [Zaghloul 2007]. The value given in Table 4 in the same reference is calculated using the asymptotic expression for $y \to 0$

$ Calculated using the asymptotic expression for $y \to 0$

Table 3: Values of the relative errors $\delta_V = |(V - V_{ref})|/V_{ref}$ & $\delta_L = |(L - L_{ref})/L_{ref}|$ in calculating the Faddeyeva function by different codes using values of the function calculated using erfi*(z)* from Mathematica™ as reference values.

| z | | Algorithm 680 | | Zaghloul [2007] | Present algorithm | | |
|---|---|---|---|---|---|---|---|
| x | y | V | L | V | V | L | Maximum no. of series terms |
| 6.3e-002 | 1.0e-020 | 0 | 2.0e-016 | 0 | 0 | 2.0e-15 | 12 |
| 6.3e-002 | 1.0e-014 | 0 | 2.0e-016 | 1.1e-016 | 1.1e-016 | 2.1e-15 | 12 |
| 6.3e-002 | 1.0e-012 | 0 | 2.0e-016 | 0 | 0 | 2.1e-15 | 12 |
| 6.3e-002 | 1.0e-010 | 1.1e-016 | 2.0e-016 | 1.1e-016 | 0 | 2.1e-15 | 12 |
| 6.3e-002 | 1.0e-006 | 0 | 0 | 2.2e-016 | 1.1e-016 | 3.7e-15 | 12 |
| 6.3e-002 | 1.0e-002 | 1.1e-016 | 2.0e-016 | 2.3e-016 | 3.4e-016 | 2.2e-15 | 12 |
| 6.3e-002 | 1.0e+001 | 2.5e-016 | 1.6e-016 | 0 | 1.2e-016 | 3.6e-15 | 13 |
| 6.3e-002 | 1.2e+001 | 3.0e-016 | 2.2e-016 | 3.0e-016 | 3.0e-016 | 8.9e-16 | 13 |
| 6.3e-002 | 1.5e+001 | 3.7e-016 | 6.9e-016 | 0 | 5.5e-016 | 2.6e-15 | 13 |
| 6.3e-002 | 2.0e+002 | 3.1e-016 | 1.2e-016 | 1.1e-015 | 0 | 2.4e-15 | 13 |
| | | | | | | | |
| 6.3e+000 | 1.0e-020 | **1.0e-000** | 3.1e-016 | **2.6e-005** | **2.4e-015** | 1.2e-016 | 13 |
| 6.3e+000 | 1.0e-014 | **3.8e-002** | 3.1e-016 | **0** | **5.3e-015** | 1.2e-016 | 13 |
| 6.3e+000 | 1.0e-012 | **3.9e-004** | 3.1e-016 | **0** | **6.2e-015** | 1.2e-016 | 13 |
| 6.3e+000 | 1.0e-010 | **3.9e-006** | 1.5e-016 | **0** | **6.1e-015** | 1.2e-016 | 13 |
| 6.3e+000 | 1.0e-006 | **3.9e-010** | 4.6e-016 | **0** | **6.0e-015** | 1.2e-016 | 13 |
| 6.3e+000 | 1.0e-002 | 6.2e-014 | 3.1e-016 | 7.3e-016 | 5.3e-015 | 2.4e-016 | 13 |
| 6.3e+000 | 1.0e+001 | 0 | 4.1e-016 | 3.4e-016 | 6.2e-015 | 2.4e-016 | 13 |
| 6.3e+000 | 1.2e+001 | 0 | 5.4e-016 | 9.4e-016 | 6.2e-015 | 0 | 13 |
| 6.3e+000 | 1.5e+001 | 2.2e-016 | 7.8e-016 | 4.3e-016 | 6.1e-015 | 1.2e-016 | 13 |
| 6.3e+000 | 2.0e+002 | 3.1e-016 | 3.1e-016 | 1.2e-007 | 6.3e-015 | 0 | 13 |
| | | | | | | | |
| 6.3e+002 | 1.0e-020 | 6.1e-016 | 0 | 3.8e-006 | 6.1e-016 | 1.2e-016 | 13 |
| 6.3e+002 | 1.0e-014 | 8.5e-016 | 0 | 3.8e-006 | 8.5e-016 | 1.2e-016 | 13 |
| 6.3e+002 | 1.0e-012 | 6.8e-016 | 0 | 3.8e-006 | 6.8e-016 | 1.2e-016 | 13 |
| 6.3e+002 | 1.0e-010 | 6.9e-016 | 0 | 3.8e-006 | 6.9e-016 | 1.2e-016 | 13 |
| 6.3e+002 | 1.0e-006 | 0 | 0 | 3.8e-006 | 5.7e-016 | 1.2e-016 | 13 |
| 6.3e+002 | 1.0e-002 | 7.0e-016 | 0 | 3.8e-006 | 7.0e-016 | 2.4e-016 | 13 |
| 6.3e+002 | 1.0e+001 | 7.2e-016 | 2.422e-016 | 3.8e-006 | 0 | 2.4e-016 | 13 |
| 6.3e+002 | 1.2e+001 | 6.0e-016 | 2.422e-016 | 3.8e-006 | 6.0e-016 | 0 | 13 |
| 6.3e+002 | 1.5e+001 | 4.8e-016 | 0 | 3.8e-006 | 4.8e-016 | 1.2e-016 | 13 |
| 6.3e+002 | 2.0e+002 | 0 | 0 | 3.0e-006 | 0 | 0 | 13 |
| | | | | | | | |
| 1.0e+000 | 1.0e-020 | 0 | 0 | 0 | 0 | 1.8e-016 | 13 |
| 5.5e+000 | 1.0e-014 | **1.2e-004** | 0 | **0** | **0** | 0 | 13 |

Table 4: Accuracy vs efficiency trade-off of the present algorithm (computations are performed on an array of 2,840,710 points generated using the grid $y=logspace(-20, 4, 71)$ and $x=linspace(-200, 200, 40001)$) using **Matlab 7.9.0.529 (R2009b)** ).

| tiny | $\delta_{V,max} = \left|V_{tiny} - V_{tiny_{min}}\right|/V_{tiny_{min}}$ | $\delta_{L,max} = \left|L_{tiny} - L_{tiny_{min}}\right|/L_{tiny_{min}}$ | Run time (s) |
|---|---|---|---|
| 0.06447×ε | 0 | 0 | 8.42 |
| 1.0e-15 | 4.6e-013 | 4.6e-013 | 8.01 |
| 1.0e-14 | 4.6e-013 | 1.2e-013 | 7.74 |
| 1.0e-12 | 2.8e-012 | 8.9e-011 | 7.10 |
| 1.0e-10 | 2.7e-010 | 6.8e-009 | 6.53 |
| 1.0e-09 | 2.6e-009 | 5.8e-008 | 6.27 |
| 1.0e-08 | 2.5e-008 | 4.9e-007 | 5.83 |
| 1.0e-07 | 2.4e-007 | 4.2e-006 | 5.56 |
| 1.0e-06 | 2.3e-006 | 3.5e-005 | 5.46 |
| 1.0e-05 | 2.2e-005 | 3.1e-004 | 5.06 |
| 1.0e-04 | 2.1e-004 | 3.1e-003 | 4.66 |

Table 5: Running times of the present algorithm (for three values of the parameter *tiny*) compared with other competitive algorithms (computations are performed on an array of 2,840,710 points generated using the grid $y=logspace(-20, 4, 71)$ and $x=linspace(-200, 200, 40001)$) using **Matlab 7.9.0.529 (R2009b)**)[*]

| Algorithm | Run time (s) | Comments |
|---|---|---|
| Faddeyeva, *tiny*=0.06447×ε | 8.42 | |
| Faddeyeva, *tiny*=1e-8 | 6.27 | |
| Faddeyeva, *tiny*=1e-4 | 4.66 | |
| Poppe & Wijers [1990] | 107.43 | ▪ Large error in the vicinity of $x$=6.3 & very small values of $y$ |
| Humlíček [1982] (original) | 23.21 | ▪ Large error and loss of accuracy in the vicinity of $x$=5.6 and very small values of $y$ |
| Humlíček [1982] (modified) | 9.87 | |
| Weidemann [1994], N=16 | 1.78 | ▪ Negative values for $V(x,y)$ near $x$-axis |
| Weidemann [1994], N=32 | 2.57 | ▪ Incorrect behavior and order of magnitude of $\partial V(x,y)/\partial x$ for very small values of $y$. |
| Weidemann [1994], N=64 | 6.88 | |
| Weidemann [1994], N=128 | 7.09 | |
| Weidemann [1994], N=256 | 12.7 | |
| Hui et al [1978] | 0.42 | ▪ Large error for small values of $y$<br>▪ Negative values for $V(x,y)$ (e.g. at $y=10^{-5}$ & $x$=4)<br>▪ Incorrect behavior and order of magnitude of $\partial V(x,y)/\partial x$ for very small values of $y$. |

∗ Timing results depend on both hardware and the version of the software used and can change significantly.

## ACKNOWLEDGMENTS

The authors would like to acknowledge valuable comments and suggestions received from the reviewers. In particular the comments and suggestions received from the associate editor, algorithm editor and the fourth anonymous referee were extremely helpful and insightful. We would like also to thank Prof. C. Benner from College of William and Mary, Williamsburg, VA, USA, for sending us a copy of Letchworth & Benner's computer code.